\documentclass[absolute]{mymathart}
\usepackage[theorems]{mymathmacros}
\usepackage[usenames,dvipsnames]{color}
\usepackage{subfig}
\usepackage[shortlabels]{enumitem}
\usepackage{hyperref}
\usepackage{amsmath}
\usepackage{amssymb}
\usepackage{graphicx}
\usepackage{epigraph}
\usepackage{amsfonts,amssymb,color}
\usepackage[mathscr]{eucal}
\usepackage{amsmath, amsthm}
\usepackage{mathrsfs}
\usepackage{amsbsy}
\usepackage{float}
\usepackage{verbatim}
\usepackage{amsfonts} 
\usepackage{pifont}

\title[A remark on Liao and Rams' result]{A remark on Liao and Rams' result on distribution of the leading partial quotient with growing speed $e^{n^{1/2}}$ in continued fractions}

\author{Liangang Ma}
\address{Dept.\ of Mathematical Sciences, University of Liverpool, Liverpool L69 7ZL, UK}
\email{maliangang000@163.com}
\newtheorem*{theorem1}{Theorem 1}
\newtheorem*{theorem2}{Theorem 2}
\newtheorem*{lemma1}{Lemma 1}
\newtheorem*{corollary1}{Corollary 1}
\newtheorem*{lemma2}{Lemma 2}
\newtheorem*{remark1}{Remark 1}

\numberwithin{equation}{section}

\begin{document} 
 \begin{abstract}
 For a real $x\in(0,1)\setminus\mathbb{Q}$, let $x=[a_1(x),a_2(x),\cdots]$ be its continued fraction expansion. Denote by
\begin{center}
$T_n(x):= max \{a_k(x): 1\leq k\leq n\}$
\end{center} 
the leading partial quotient up to $n$. For any real $\alpha\in(0,\infty), \gamma\in(0,\infty)$, let
\begin{center}
$F(\gamma,\alpha):=\{x\in(0,1)\setminus\mathbb{Q}: \lim_{n\rightarrow\infty}\frac{T_n(x)}{e^{n^\gamma}}=\alpha\}$.
\end{center}
For a set $E\subset (0,1)\setminus\mathbb{Q}$, let $dim_H E$ be its Hausdorff dimension. Recently Lingmin Liao and Michal Rams \cite[Theorem 1.3]{LR} show that 
\begin{center}
$ dim_H F(\gamma,\alpha)=\left\{
\begin{array}{ll}
1 &\ \ if\ r\in(0,1/2) \\
1/2 &\ \  if\ r\in(1/2,\infty) \\
\end{array}
\right.$ 
\end{center}
for any $\alpha\in(0,\infty)$. In this paper we show that $dim_H F(1/2,\alpha)=1/2$ for any $\alpha\in(0,\infty)$ following Liao and Rams' method, which supplements their result.

 \end{abstract}
 
 \maketitle

Through out the paper we follow Liao and Rams' notations \cite{LR}.  As mentioned in the abstract, we aim to show that 

\begin{theorem1}
$dim_H F(1/2,\alpha)=1/2$.
\end{theorem1}
We only prove $dim_H F(1/2,1)=1/2$, as one can show the theorem for any $\alpha\in\mathbb{R}^+:=(0,\infty)$ by the same process. In order to do this, we first show that
\begin{lemma1}
Let $L\in\mathbb{R}^+$ be a constant. Let $n_k:=[(\frac{k}{L})^2]$ (the integer part of $(\frac{k}{L})^2$), $k\in\mathbb{N}$. Then for any $x\in F(1/2,1)$ and $k$ large enough, there exists an integer $j_k, n_{k-1}< j_k\leq n_k$, such that
\begin{center}
$T_{n_k}(x)=a_{j_k}(x)$.
\end{center}
\end{lemma1}
\begin{proof}
We prove this by reduction to absurdity. Suppose there exist infinitely many integers $k_i, j_{k_i}, i\in\mathbb{N}, k_i>k_{i-1}, j_{k_i}\leq n_{k_i-1}$, such that
\begin{center}
$T_{n_{k_i}}(x)=a_{j_{k_i}}(x)$
\end{center}
for some $x\in F(1/2,1)$. Note that in this case we have
\begin{center}
$T_{n_{k_i-1}}(x)=a_{j_{k_i}}(x)$.
\end{center} 
Then for the sequence $\{n_{k_1-1}, n_{k_2-1},\cdots\}$, we have
\begin{center}
$\lim_{i\rightarrow\infty}\cfrac{T_{n_{k_i-1}}(x)}{e^{n_{k_i-1}^{1/2}}}=\lim_{i\rightarrow\infty}\cfrac{T_{n_{k_i}}(x)}{e^{[(k_i-1)^2/L^2]^{1/2}}}=\lim_{i\rightarrow\infty}\cfrac{T_{n_{k_i}}(x)}{e^{n_{k_i}^{1/2}}}\cfrac{e^{[k_i^2/L^2]^{1/2}}}{e^{[(k_i-1)^2/L^2]^{1/2}}}=1\cdot e^{1/L}\neq 1$
\end{center}
which contradicts the fact that 
\begin{center}
$\lim_{k\rightarrow\infty}\frac{T_k(x)}{e^{k^{1/2}}}=1$
\end{center}
as $x\in F(1/2,1)$. So our conclusion holds for any sufficiently large $k$.
\end{proof}
In the following we will omit the integer notation $[\ ]$ for simplicity as results will not be affected. By this lemma,
\begin{corollary1}\label{corollary1}
For $x\in F(1/2,1)$ and $n_k:=(\frac{k}{L})^2$, we have
\begin{center}
$(1-\epsilon)e^{k/L}\leq S_{n_k}(x)-S_{n_{k-1}}(x)\leq (1+\epsilon)(\frac{k}{L})^2e^{k/L}$
\end{center}
for a small $\epsilon\in\mathbb{R}^+$ and any $k$ large enough.
\end{corollary1} 

The rest of the work goes the same process as estimation of the upper bound for $E_{\varphi}$ when $\gamma=1/2$ in \cite[Proof of Theorem 1.1]{LR}. For the length of the rank-$n$ fundamental interval 
\begin{center}
$I_n(a_1,\cdots,a_n):=\{x\in(0,1)\setminus\mathbb{Q}: a_1(x)=a_1,\cdots,a_n(x)=a_n\}$,
\end{center}
we have
\begin{center}
$\prod_{i=1}^n\frac{1}{(a_i+1)^2} I_n(a_1,\cdots,a_n)\leq \prod_{i=1}^n\frac{1}{a_i^2}$.
\end{center}
Let 
\begin{center}
$A(m,n):=\{(i_1,\cdots,i_n)\in\{1,\cdots,m\}^n: \sum_{j=1}^n i_j=m\}$.
\end{center}
Let $\zeta(\cdot)$ be the Riemann zeta function. Now we quote \cite[Lemma 2.1]{LR} as following. 
\begin{lemma2}
For $s\in(1/2,1)$ and $m\geq n$, we have
\begin{center}
$\sum_{(i_1,\cdots,i_n)\in A(m,n)}\prod_{j=1}^n\cfrac{1}{i_j^{2s}}\leq (\cfrac{9}{2}(2+\zeta(2s)))^n\cfrac{1}{m^{2s}}$.
\end{center}

\end{lemma2}
Now we are in a position to bound Hausdorff dimension of $F(1/2,1)$ above.
\begin{theorem2}
$dim_H F(1/2,1)\leq 1/2$.
\end{theorem2}
\begin{proof}
Let $D_l$ be the integers in the interval $[(1-\epsilon)e^{l/L},(1+\epsilon)(\frac{l}{L})^2e^{l/L}]$. Let 
\begin{center}
$B(1/2,N):=\{\cup_{k=N}^\infty I_{n_k}(a_1,a_2,\cdots,a_{n_k}): \sum_{j=n_{l-1}+1}^{n_l} a_j=m\ \mbox{with}\ m\in D_l, N\leq l\leq k\}.$
\end{center}
By Corollary 1 one can see that
\begin{center}
$F(1/2,1)\subset\cup_{N=1}^\infty B(1/2,N)$.
\end{center}
Now we show that $dim_H B(1/2,1)\leq 1/2$. Similar method implies $dim_H B(1/2,N)\leq 1/2$ for any $N\in\mathbb{N}$, which is enough to prove our Theorem 2. By Lemma 2,
\begin{center}
$\sum_{I_{n_k}\subset B(1/2,1)}|I_{n_k}|^s\leq \prod_{l=1}^k\sum_{m\in D_l}(\cfrac{9}{2}(2+\zeta(2s)))^{n_l-n_{l-1}}\cfrac{1}{m^{2s}}$.
\end{center}
Note that $|D_l|\leq (1+\epsilon)(\frac{k}{L})^2e^{k/L}, m>(1-\epsilon)e^{k/L}$, so
\begin{center}
$
\begin{array}{lll}
 & \sum_{I_{n_k}\subset B(1/2,1)}|I_{n_k}|^s\\
\leq & \prod_{l=1}^k (1+\epsilon)(1-\epsilon)^{2s}(l/L)^2e^{(1-2s)l/L}(\cfrac{9}{2}(2+\zeta(2s)))^{\frac{2l-1}{L^2}}\\
\leq & \prod_{l=1}^k \Big(\big((1+\epsilon)(1-\epsilon)^{2s}(l/L)^2\big)^{1/l} e^{(1-2s)/L}(\cfrac{9}{2}(2+\zeta(2s)))^{3/L^2}\Big)^l.
\end{array}
$
\end{center}
Solve the equation
\begin{center}
$\cfrac{9}{2}(2+\zeta(2s))=\cfrac{1}{2}e^{\frac{2s-1}{3}L}$
\end{center}
regarding the main terms, we get a unique solution $s_L\in(1/2,1)$ when $L$ is large enough. $s_L\rightarrow 1/2$ as $L\rightarrow\infty$ since $\zeta(2\cdot\frac{1}{2})=\zeta(1)=\infty$. Then $\sum_{I_{n_k}\subset B(1/2,1)}|I_{n_k}|^s<\infty$, which forces $dim_H B(1/2,1)\leq 1/2$.

\end{proof}
As $dim_H F(1/2,\alpha)\geq 1/2$ (see \cite[Proof of Theorem 1.3]{LR}), Theorem 1 follows directly from Theorem 2.

\begin{remark1}
Our corollary 1 sharpens estimation on $S_{n_k}(x)-S_{n_{k-1}}(x)$ in \cite[Proof of Theorem 1.3]{LR} for $x\in F(1/2,1)$. In fact we can do similar estimations for any $x\in F(\gamma, \alpha), \gamma\in(0,\infty),\alpha\in\mathbb{R}^+, n_k=k^{1/\gamma}$. This enables us to give better estimation on $\sum_{I_{n_k}\subset B(\gamma,N)}|I_{n_k}|^s,$ $\gamma\in [1/2,1)$. By virtue of it, when estamating the upper bound in \cite[Proof of Theorem 1.3]{LR} for H-dimension of $F(\gamma,\alpha), \gamma\in (1/2,1)$, we can simply take $n_k=k^{1/\gamma}$ instead of $k^{1/\gamma}(\log k)^{1/\gamma^2}$.    
\end{remark1}


\begin{thebibliography}{88}


\bibitem[LR]{LR} Lingmin Liao and Michal Rams, Subexponentially increasing sums of partial
quotients in continued fraction expansions, Math. Proc. Cambridge Philos. Soc. 160, no. 03, 401-412, 2016.
\end{thebibliography}
\end{document}